\documentclass[reqno]{amsart}

\usepackage[dvipsnames]{xcolor}
\usepackage{graphicx}
\usepackage{epsfig}
\usepackage{tikz}
\usepackage{enumerate}
\usepackage{todonotes}

\usepackage[utf8]{inputenc}

\usepackage{latexsym}
\usepackage{amssymb}
\usepackage{amsmath}
\usepackage{amsthm}
\usepackage{amscd}
\usepackage[normalem]{ulem} 

\theoremstyle{plain}
\newtheorem{theorem}{Theorem}[section]
\newtheorem*{thm*}{Theorem}
\newtheorem{lemma}[theorem]{Lemma}
\newtheorem{definition}[theorem]{Definition}
\newtheorem{proposition}[theorem]{Proposition}
\newtheorem{prop*}{Proposition}

\theoremstyle{definition}
\newtheorem{remark}[theorem]{Remark}

\newtheorem{question}[theorem]{Question}
\newtheorem{corollary}[theorem]{Corollary}

\makeatletter
\def\blx@maxline{77}
\makeatother

\newcommand{\R}{\mathbb{R}} 
\newcommand{\Z}{\mathbb{Z}} 
\newcommand{\Prob}{\mathbb{P}} 

\newcommand{\e}{\epsilon}
\newcommand{\w}{\omega}
\newcommand{\W}{\Omega}

\newcommand{\AND}{\textrm{ and }}

\newcommand{\E}{\mathbb{E}}

\newcommand{\specialpts}{\mathcal{S}} 
\newcommand{\biinfinite}{bi-infinite}
\newcommand{\Biinfinite}{Bi-infinite}
\newcommand{\bigeodesics}{bigeodesics}
\newcommand{\directions}{\operatorname{\mathbf{A}}}

\newcommand{\origin}{\mathbf{0}}

\newcommand{\tmap}{T}

\def\walk{X}
\newcommand{\walks}{\mathcal{W}}
\newcommand{\rect}{\operatorname{Rect}}
\newcommand{\ver}{\operatorname{Vert}}

\DeclareMathOperator{\argmin}{argmin}

\newcommand{\Eqref}[1]{\eqref{#1}}
\newcommand{\Figref}[1]{Fig.~\ref{#1}}

\newcommand{\Secref}[1]{Section~\ref{#1}}
\newcommand{\Thmref}[1]{Theorem~\ref{#1}}
\newcommand{\Corref}[1]{Corollary~\ref{#1}}
\newcommand{\Lemref}[1]{Lemma~\ref{#1}}
\newcommand{\Propref}[1]{Prop.~\ref{#1}}
\newcommand{\Defref}[1]{Definition~\ref{#1}}
\newcommand{\Remref}[1]{Remark~\ref{#1}}

\usepackage{versions}
\includeversion{note}       
\markversion{note}
\markversion{arjunnote}
\markversion{oldnote}
\markversion{longcalc}

\excludeversion{openquestion}
\excludeversion{longcalc}
\excludeversion{note}
\excludeversion{arjunnote}
\excludeversion{oldnote}



\usepackage[nonotes]{optional}

\newcommand{\arjunhl}[1]{\opt{finalnotes}{\textcolor{blue}{#1}}} 
\newcommand{\arjunnotes}[1]{\opt{arjunnotes}{\textcolor{OliveGreen}{#1}}}

\newcommand{\oldnotes}[1]{\opt{oldnotes}{\textcolor{orange}{#1}}}

\title{Stationary coalescing walks on the lattice}
\author{Jon Chaika \and Arjun Krishnan}
\keywords{measure preserving transformations, bi-infinite trajectories, geodesics, first-passage percolation}

\date{}

\usepackage[sorting=none,backend=biber,natbib=true,style=numeric]{biblatex}
\usepackage[urlcolor=blue,colorlinks=true,linkcolor=blue,citecolor=blue]{hyperref}
\begin{document}

\begin{abstract}
We consider translation invariant measures on families of nearest-neighbor semi-infinite walks on the integer lattice. We assume that once walks meet, they coalesce. In $2d$, we classify the collective behavior of these walks under mild assumptions: they either coalesce almost surely or form bi-infinite trajectories. Bi-infinite trajectories form measure-preserving dynamical systems, have a common asymptotic direction in $2d$, and possess other nice properties. We use our theory to classify the behavior of compatible families of semi-infinite geodesics in stationary first- and last-passage percolation. We also partially answer a question raised by C.\ Hoffman about the limiting empirical measure of weights seen by geodesics. We construct several examples: our main example is a standard first-passage percolation model where geodesics coalesce almost surely, but have no asymptotic direction or average weight.
\end{abstract}

\maketitle
\tableofcontents

\section{Introduction}
\label{sec:intro}
Let $(\W,\mathcal{F},\Prob,\{T^z\}_{z \in \Z^d})$ be a $\Z^d$ measure-preserving dynamical system. Let $\walks(\w)$ be a stationary or translation-covariant subset ---that is, $\walks(T^z\w) = \walks(\w) - z$--- of the lattice $\Z^d$, and suppose that it contains the origin with positive probability; i.e., $\Prob( \origin \in \walks(\w)) > 0$. Consider a family of measurable walks on the lattice $\{\walk_z\}_{z \in \walks}$, where each $\walk_z \colon \Omega \times \Z^+ \to \Z^d$ is a nearest-neighbor path that starts at $z$. We assume that almost surely for all $k \in \Z^+$ and $z \in \walks(\w)$, these walks have been created in a stationary way: 
\[
    \walk_z(\w, k) = x + \walk_{z-x}(\tmap^x \w, k),
\]
\arjunnotes{Modified the above display. It originally read $\walk_z(T^x\w, k) = \walk_{x+z}(\w, k)$}
and that they are \emph{compatible}: 
\begin{align*}
    X_z(\w,k+1) = X_{X_z(\w,1)}(\w,k) .
\end{align*}
The compatibility condition implies that if two walks meet at a point at some time, then they remain together in the future; i.e., the two walks must coalesce and \emph{cannot cross each other}. 

Because walks coalesce when they meet, we may assume that there is a stationary vector-field $\alpha$ that is the discrete time-derivative of the walks:
\begin{equation}
    \alpha(\w,z) = \walk_z(\w,1) - \walk_z(\w,0).
    \label{eq:direction function definition}
\end{equation}
The $\alpha$ function takes values in $\directions \subset \{ \pm e_1, \ldots \pm e_d \}$, and we call a particular $\alpha$ value an \emph{arrow}. We study walks that do not form loops, and satisfy a mild line-crossing assumption (see \Thmref{thm:dichot finish}); these conditions are automatically satisfied by \emph{directed} walks, where $\directions = \{e_1,\ldots,e_d\}$. One could think of the walks as the flow generated by the stationary vector field of arrows. We call an arrow configuration \emph{non-trivial} if $\alpha$ is not constant almost surely. The canonical walk $\walk(\w)$ starts at the origin and $\alpha(\w)$ is its (discrete) derivative at time $0$. We will omit the $\w$ from the notation when it is clear from context. 

We frequently speak of \emph{configurations} on the lattice: for any $\w$, this refers to the collection of walks $\{\walk_z(\w)\}_{z \in \walks}$ or equivalently, the collection of arrows $\{ \alpha(T^z \w) \}_{z \in \walks}$. 

\begin{figure}[!ht]
    \label{fig:periodic example}
\begin{center}
    \includegraphics[width=3cm]{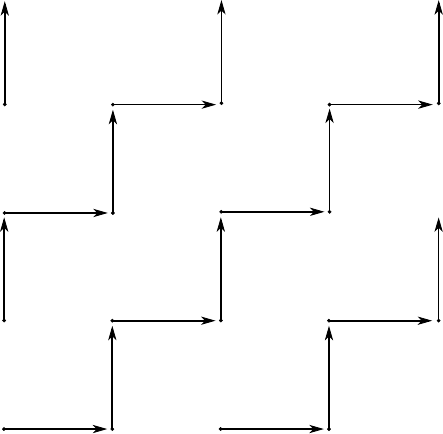}
\end{center}
\caption{The space is $\W = \{\w_1,\w_2\}$ with uniform measure. The arrows are given by $\alpha(\w_i) = e_i$ for $i = 1,2$. The translation operators $T^{e_i}$ simply swap between $\w_1$ and $\w_2$.}
\end{figure}

As a first example, consider independent and identically distributed (iid) arrows taking values in $\directions = \{e_1,e_2\}$ (with probabilities $p$ and $1-p$) on each point of the lattice $\Z^2$. Each walk is a classical simple random walk, where $e_1$ corresponds to stepping up and $e_2$ corresponds to stepping down. Consider any two random walks starting at $x \neq y$ on an anti-diagonal of the form $\{ z \colon z_1 + z_2 = c\}$ for some fixed constant $c$. As long as $X_x(\omega,k) \neq X_y (\omega,k)$, the projection of the difference $\left(\walk_x(\w,k) - \walk_y(\w,k)\right)\cdot(-1,1)$ is a one-dimensional simple random walk where time proceeds along the main diagonal in $\Z^2$. The $k$\textsuperscript{th} step of the walk involves arrows on the antidiagonal line $\{ z \colon z_1 + z_2 = c + k-1\}$. Unless the two walks have coalesced previously, the $k$\textsuperscript{th} step is independent of the previous steps, takes the values $\pm2$ with equal probability $p(1-p)$, and is $0$ otherwise. Hence $\left( \walk_x - \walk_y \right)\cdot(-1,1)$ is almost surely recurrent to $0$, and the walks must coalesce. Thus, every pair of walks from points $x,y \in \Z^2$ must coalesce almost surely. In contrast, the periodic system in~\Figref{fig:periodic example} (also a measure-preserving ergodic $\Z^2$ system) has \biinfinite{} (see \Defref{def:biinfinite trajectory definition}) trajectories that do not coalesce.

In this paper, we completely classify the collective behavior of the trajectories in $d=2$ under mild assumptions. There is a behavioural \emph{dichotomy} (Theorem \ref{thm:dichot finish}): with probability $1$,
\begin{enumerate}
    \item the walks from all $x,y \in \Z^2$ coalesce, and no \biinfinite{} trajectories exist, or 
    \item a positive fraction of the walks in a configuration form \biinfinite{} trajectories that are themselves measure-preserving dynamical systems. Thus, all the walks have the same asymptotic direction (Theorem \ref{thm:measure preserving dynamical system and asym dir}), and no two \biinfinite{} trajectories coalesce (Theorem \ref{thm:no spec coal}). In fact, all other trajectories must eventually coalesce with \biinfinite{} trajectories (Corollary \ref{cor:walks that are not part of biinf}).
\end{enumerate}

There is a sequel to this paper that explores various entropic properties of \biinfinite{} trajectories. For example, we show that in (factors of) iid systems, \biinfinite{} trajectories must carry entropy. We also construct a discrete symmetric simple exclusion process that has \biinfinite{} trajectories carrying entropy.

All the nice properties that the \biinfinite{} trajectories possess are not shared by almost surely coalescing walks. For example, asymptotic velocity is no longer guaranteed. We demonstrate this by constructing an explicit example (\Thmref{thm:coal but no asym direction in z2}).

In dimensions higher than $2$, the \biinfinite{} trajectories/almost-sure coalescence dichotomy is not true. We construct an example (Corollary \ref{thm:both almost sure coal and biinf walks}) in $d = 3$ where almost surely,

\begin{enumerate}[(i)]
    \item every trajectory does not have an asymptotic direction,
    \item every configuration does not have bi-infinite trajectories, and
    \item we do not have almost sure coalescence. 
\end{enumerate}

\subsection{First- and last-passage percolation}
Our model is motivated by questions about the behavior of infinite geodesics in first- and last-passage percolation. Let $\W = \{\w_z \in \R\}_{z \in \Z^d}$ with product $\sigma$-algebra and a translation invariant measure $\Prob$. The $\w_z$ are called \emph{weights} and they are typically nonnegative random variables. Let $\walk_{x,y}$ be a path from $x$ to $y$ and let the total weight of the path be the sum $W(\walk_{x,y}) := \sum_{z \in \walk_{x,y}} \w_z$. Define the first-passage time from $x$ to $y$ to be
\[
    T(x,y) = \inf_{\walk_{x,y}} W(\walk_{x,y}).
\]
{The models may have weights on the edges of $\Z^d$ instead of the vertices.} The first-passage time $T(x,y)$ satisfies a triangle inequality; if the weights are strictly positive, it defines a random metric on the lattice $\Z^d$. A geodesic for this random metric is a nearest-neighbor path that minimizes the passage time between every vertex that lies on it.

By considering the geodesic from $0$ to $n e_1$, and then looking at a subsequence as $n \to \infty$, it is clear that there is at least one semi-infinite geodesic from the origin. Furstenberg (communicated in~\citet{kesten_aspects_1986} page 134) asked if there exist \biinfinite{} geodesics (bigeodesics) in first-passage percolation with iid weights. This question has not been answered completely. In $d \geq 3$, very little progress has been made on geodesic behavior, but for $d = 2$, however, there are several partial answers under different assumptions on the so-called time-constant. We survey a few important results below. 

For $u \in S^1$, define the time-constant $g(u)$ of first-passage percolation as
\[
    g(u) = \lim_{n \to \infty} \frac{T(0,[nu])}{n},
\]
where $[nu]$ represents the closest lattice point to $nu$, with some chosen way of breaking ties. Assuming that the weights are $L^1$, the limit exists almost surely and in $L^1$ by the  subadditive ergodic theorem \citep{kingman_ergodic_1968}. 

In a seminal paper~\citep{MR1387641}, Licea and Newman prove theorems about the non-existence of \bigeodesics{} under strong assumptions on the time-constant. They assume that 1) the level sets of $g(u)$ are ``uniformly curved'', 2) the weights satisfy a property called finite-energy (see~\citet{MR990777}) and 3) the weights have continuous distribution. Then, with probability $1$, except for a deterministic Lebesgue measure $0$ set of directions $u \in S^1$, 
\begin{enumerate}
    \item there exists exactly one geodesic from each point $x \in \Z^2$ in direction $u$,
    \item geodesics from different points coalesce almost surely, and
    \item there are no \bigeodesics{} in direction $(-u,u)$.
\end{enumerate}

Busemann functions are a useful tool in the analysis of infinite geodesics. It is defined as the limit (if it exists)
\begin{equation}
    B_u(x,y) = \lim_{n \to \infty} T(x,z_n) - T(y,z_n)
\end{equation}
such that $z_n/n \to u \in S^1$. \citet{newman_surface_1995} showed the existence of the above limit under certain strong hypotheses. Notably, Hoffman \citep{MR2114988} realized their importance, and used Busemann functions to prove (concurrently with~\citep{MR2115045}) that there are at least two semi-infinite geodesics under no assumptions on the time-constant $g(u)$. Busemann functions have many useful properties, but from our perspective, the most useful property is that they \emph{encode geodesic behavior in a stationary manner}. Moreover, Busemann geodesics coalesce when they meet. These two properties motivate our assumptions about the stationary coalescing walks we consider.

In first- and last-passage percolation, dual variational descriptions of the time-constant $g(u)$ have recently been proved by~\citet{krishnan_cpam_variational_2016} and \citet{MR3535900}. Here, the time constant is expressed as a minimization problem over \emph{functions} instead of paths. It turns out that certain special minimizers of the formula are the Busemann functions. In stochastic homogenization, Busemann functions are known as correctors.

Damron and Hanson built a theory of generalized Busemann functions in first-passage percolation~\citep{MR3152744}. \arjunnotes{Theirs and previous work are summarized in the survey by \citet{auffinger_50_2015}.}Associated with each direction $u \in S^1$, they construct a stationary function $B_u(x,y) \colon \Z^2 \times \Z^2 \to \R$ called a reconstructed Busemann function. By building geodesics associated with these reconstructed Busemann functions, they obtain a directed geodesic graph $G_u$ with vertices in $\Z^2$ such that
\begin{enumerate}
    \item each directed path on $G_u$ is a geodesic,
    \item if there is a path from $x$ to $y \in \Z^2$ in $G_u$, then $B_u(x,y) = T(x,y)$,
    \item $G_u$ has no loops even as an undirected graph,
    \item each vertex has out degree $1$.
\end{enumerate}
The edges of this directed graph are the analogs of our arrow configurations $\{\alpha(T^z \w)\}$. They then show that under the upward finite-energy assumption \citep{MR3152744}, the geodesic graph coalesces almost surely. \citet{2016arXiv160902447A} prove similar results using different methods under weaker assumptions. Our dichotomy theorem classifies the behavior of such geodesic walks even when finite-energy does not hold (see \Corref{cor:dichotomy theorem for first and last passage percolation graphs}). 

\citet[Theorem 2.3]{georgiou_geodesics_2015} prove the last-passage version of the results of~\citep{MR3152744} in $d=2$: If the last-passage time-constant $g$ is differentiable at $u \in S^1$, there exists a stationary function $B_u(\w,x,y)$ that can be used to create a family of compatible semi-infinite geodesics going in direction $u$ by following the arrows defined by
\begin{equation}
    \alpha(x,\w) = \argmin_{z = e_1,e_2}B_u(\w,x,x + z),
    \label{eq:arrows from min of busemann function}
\end{equation}
with some specified way of breaking ties in the $\argmin$. Again, under the assumption of finite-energy, they prove almost sure coalescence of geodesics using the Licea-Newman argument. When the underlying graph is changed, finite-energy is no longer enough to prove almost sure coalescence. \citet{benjamini_first_2016} show that first-passage percolation on hyperbolic graphs with iid weights have \bigeodesics{}.

The property that one can create geodesics using~\Eqref{eq:arrows from min of busemann function} is rather special. In the general setting with stationary-ergodic weights, minimizers of the variational formulas with this property do not necessarily exist\footnote{see~\citet{lions_correctors_2003} in the context of continuum stochastic homogenization}. Nevertheless, one can still create paths (that are not necessarily geodesics) from any minimizer using the recipe in~\Eqref{eq:arrows from min of busemann function}. Our paper is an attempt to classify the behavior of these paths to better understand the behavior of the variational formulas of first- and last-passage percolation.

A different set of results that follow from our theorems comes from a question asked by C.\,Hoffman at the American Institute of Mathematics workshop in 2016: ``On a semi-infinite geodesic, does the empirical measure of weights seen on the geodesic converge?'' We interpret this as follows. Given a family of compatible geodesics $\{\walk_x\}_{x \in \Z^d}$ and vertex weights $\{\w_z\}_{z \in \Z^d}$, does the measure defined by $n^{-1} \sum_{i=1}^n \delta_{\w(\walk_x(i))}$ converge as $n \to \infty$? Under an integrability condition, one can also ask if there is a limiting asymptotic weight on a geodesic; i.e., does $n^{-1} \sum_{i=1}^n \w(\walk(i))$ converge? When bigeodesics exist (\Thmref{thm:measure preserving dynamical system and asym dir}), these limits exist (the latter exists when $\w(0)$ in integrable), though the limit may not be deterministic. 

However, these limits \emph{do not exist on geodesics in general}. We use the example constructed in \Thmref{thm:coal but no asym direction in z2} to build a standard first-passage percolation model with edge-weights such that with probability $1$, there exists a positive density of points on the lattice with compatible geodesics that do not have an asymptotic direction and do not carry an empirical measure of weights. 

\subsection{Acknowledgments}
The authors thank Firas Rassoul-Agha, Timo Sepp\"al\"ainen, Eric Cator, and Michael Damron for helpful conversations. J.\ Chaika was supported in part by NSF grants DMS-135500 and DMS-1452762, the Sloan foundation and a Warnock chair. A.\ Krishnan was supported in part by an AMS Simons travel grant.

\section{Main results}
\label{sec:main results}
\noindent \textbf{Notation}
We will generally use Greek or calligraphic letters to denote events $\mathcal{E} \in \mathcal{F}$. We write subsets of $\Z^d$, random or deterministic, using the Latin alphabet. We will frequently speak of ``configurations having a density of points with a certain property''. We explain what we mean by this in the following.

\begin{definition}[Rectangular subsets of $\Z^d$]
    Let the rectangle centered at $x \in \Z^d$ with side lengths $(N_1,\ldots,N_d)$ be
    \[
        \rect_x(N_1,\ldots,N_d) = \prod_{i=1}^d [x_i - N_i, x_i + N_i] .
    \]
\end{definition}
If the side-lengths are equal, then we write $\rect_x(N)$. The boundary of any $R \subset \Z^d$ is written as $\partial R$ and consists of the set of points in $R$ that have at least one point in $\Z^d \setminus R$ as a nearest neighbor.

A subset $A \subset \Z^d$ has density if there is a number $c \in [0,1]$ such that $|\Z^d \cap \rect_0(N)|^{-1} |A \cap \rect_0(N)| \to c$. Given an event $\mathcal{E}$, the ergodic theorem ensures that almost surely in every configuration, the random subset $A(\w) = \{x \in \Z^d \colon T^x\w \in \mathcal{E} \}$ occurs with density $\Prob(\mathcal{E})$. Motivated by this interpretation, we will sometimes abuse notation and write for a set $A \subset \Z^d$ and an event $\mathcal{E}$,
\begin{equation}
    (A\cap \mathcal{E})(\w) := \{ x \in A \colon T^x \w \in \mathcal{E} \}.
    \label{eq:intersection of a set and an event}
\end{equation}

\begin{definition}[Coalescence of points]
    Given a configuration, we say that the points $x$ and $y$ coalesce if the walks $\walk_x$ and $\walk_y$ coalesce in the future. That is, for some $k_0,k_1 \in \Z^+$, $X_x(\w,k_0) = X_y(\w,k_1)$. We say we have \emph{almost sure coalescence} if almost surely for all $x,y \in \walks$, the walks through $x$ and $y$ coalesce. 
\end{definition}

\begin{definition}[\Biinfinite{} walks and points]
    Fix a point $z \in \Z^d$. We say that the point $z$ is \biinfinite{} if for every $n \in \Z^+$, there is a sequence of points $\{a_n\}_{n=0}^\infty \in \Z^d$ such that for each $n$, $\walk_{a_n}(\w,i) = a_{n-i}$ for $i=1,\ldots,n-1$ and $\walk_{a_n}(\w,n) = z$. We call this union of (one-sided) walks $\cup_{n \in \Z^+} \cup_{i=0}^{\infty} \walk_{a_n}(\w,i)$ a \biinfinite{} trajectory. 
    \label{def:biinfinite trajectory definition}
\end{definition}
\arjunnotes{This is one small thing we would have to change if wanted to move away from nearest-neighbor walks.}
\begin{remark}
   Using the ergodic decomposition, we may restrict our attention to each ergodic component. So without loss of generality, we consider the measure to be ergodic in Theorems \ref{thm:dichot finish}, \ref{cor:walks that are not part of biinf}, \ref{thm:no spec coal}, \ref{thm:measure preserving dynamical system and asym dir}, and \ref{cor:biinfinite trajectories in dim 2}.
\end{remark}
\begin{theorem}\label{thm:dichot finish}In $\mathbb{Z}^2$, suppose we have a positive density of walks $\walks$ that each
    \begin{enumerate}
        \item have no loops, and
        \item cross every vertical line a last time, and stays strictly to the right of it thereafter. Precisely, for each $z \in \walks$ and $a \in \Z$ such that $a > z \cdot e_1$, there is a $k_0$ (which depends on $a \AND z$) such that $\walk_z(\w,k_0) \cdot e_1 = a$ and $\walk_z(\w,k)\cdot e_1 > a$ for every $k > k_0$. 
    \end{enumerate} 
    Then, if we do not have almost sure coalescence, there is {a positive density of} \biinfinite{} trajectories in every configuration with probability $1$.
\end{theorem}
\begin{remark}
    Vertical lines may clearly be replaced with horizontal lines or diagonal lines $(x = \pm y)$. The crossing can be changed from right to left. We also believe that the proofs ought to generalize to any set of parallel lines, but do not pursue this here. 
    Note that if our set of arrows excludes one of $\{e_1,-e_1,e_2,-e_2\}$ then assumption (2) is automatically satisfied {with either vertical or horizontal lines}. 
\end{remark}
\arjunnotes{What prevents one walk from crossing all vertical lines and another walk from crossing all horizontal lines? Is it the ergodic theorem? Yes.}
Suppose \biinfinite{} trajectories exist with positive probability. A priori, not all walks in $\mathcal{W}$ have to be contained in some \biinfinite{} trajectory, and so we can ask about the behavior of these other walks. 

\begin{corollary}[of \Thmref{thm:dichot finish}]
    In $d=2$, suppose the walks in $\mathcal{W}$ satisfy the conditions in Theorem \ref{thm:dichot finish}, and suppose \biinfinite{} trajectories exist with positive probability. Then, walks in $\walks$ that are not part of \biinfinite{} trajectories must eventually coalesce with \biinfinite{} trajectories. 
    \label{cor:walks that are not part of biinf}
\end{corollary}
\begin{remark}
    \Corref{cor:walks that are not part of biinf} is used to show that ergodic averages converge on all trajectories in $\walks$ when \biinfinite{} trajectories exist in \Corref{cor:biinfinite trajectories in dim 2}.
\end{remark}
Next, we investigate the behavior of \biinfinite{} trajectories in some detail. The following theorem is the converse of~\Thmref{thm:dichot finish}. 
\begin{theorem}\label{thm:no spec coal} 
    Almost surely in each configuration, no two \biinfinite{} trajectories in $\walks$ may coalesce, and thus we cannot have almost sure coalescence.
\end{theorem}
    \begin{question}
        Is there a natural measure of randomness that is weaker than finite-energy (like strong-mixing or total-ergodicity) that distinguishes between the \biinfinite{} trajectories and almost sure coalescence situations? 
    \end{question}
In first- and last-passage percolation, \Thmref{thm:no spec coal} is well-known in a slightly different setting. We provide Theorem \ref{thm:no spec coal}, whose proof is analogous to the proof in that setting because it is a technical step in some of our results. It also motivates the term bi-infinite trajectory in \Defref{def:biinfinite trajectory definition} (as opposed to a graph or tree of \biinfinite{} trajectories). 
 We briefly recall the argument in first/last passage percolation. First, one constructs families of compatible semi-infinite geodesics in first- or last-passage percolation. 
  This requires unproven but reasonable assumptions on the differentiability of the time-constant. Assuming the weights satisfy the finite-energy condition, one shows that these geodesics coalesce almost surely. Second, one assumes that \biinfinite{} trajectories exist in the presence of the almost sure coalescence of these geodesics, and shows a contradiction using a Burton-Keane lack-of-space argument. Thus, for example, Theorem $4.6$ in~\cite{georgiou_geodesics_2015} and Theorem $6.9$ in~\cite{MR3152744} prove that \biinfinite{} geodesics cannot exist in these families of almost surely coalescing geodesics. However, the heart of the matter in the second part is that \biinfinite{} trajectories cannot coalesce; this is the content of \Thmref{thm:no spec coal} and completes the dichotomy. 

Next, we state a simple corollary of Theorem~\ref{thm:dichot finish} in first-passage percolation. It completely classifies the behavior of all compatible geodesic families that cross-lines in the sense of~\Thmref{thm:dichot finish}. \citet{MR3152744} constructed such families under two sets of assumptions, $\textbf{A1}$ and $\textbf{A2}$. $\textbf{A1}$ assumes iid weights and hence they automatically satisfy the finite-energy condition. So the argument of Licea and Newman shows that these geodesics coalesce almost surely. Assumption $\textbf{A2}$ does not imply finite-energy in general, and hence we do not always have almost sure coalescence.
\arjunnotes{
    \newline I have included A1 here for us.
\noindent\textbf{A1.}
\begin{enumerate}
    \item $\W = \{\w_e\}_{e \in E}$ is the set of non-negative edge-weights on $E$, the nearest-neighbor edges on $\Z^2$. $\Prob$ is the product measure.
    \item The weights satisfy the Cox-Durrett condition~\citep{cox_limit_1981}:
        \[
            \E \left[ \left( \min_{e \in \{\pm e_1,\pm e_2\}} \w_e \right)^2 \right] < \infty.
        \]
    \item $\Prob(\w_e = 0) < p_c = \frac{1}{2}$, the critical probability for bond-percolation on $\Z^2$ (c.f. \citep[Theorem 1.15]{kesten_aspects_1986}).
    \item The weights have continuous distribution. 
\end{enumerate}
}
We recall assumption \textbf{A2} below: 
\begin{enumerate}
    \item $(\W,\mathcal{F},\Prob)$ is an ergodic $\Z^2$ system of weights.
    \item $\Prob$ has all the symmetries of $\Z^2$.
    \item $\E[\w_e^{2+\e}] < \infty$ for some $\e > 0$.
    \item the limit shape for $\Prob$ is bounded.
    \item $\Prob$ has unique passage times; i.e., if $X_1$ and $X_2$ are two finite paths, $W(X_1) \neq W(X_2)$ almost surely.
\end{enumerate}
\arjunnotes{Under assumption \textbf{A1} or \textbf{A2} along with the finite-energy assumption, they prove that the compatible geodesic family they construct must coalesce almost surely~\citep[][Theorem 1.10]{MR3152744}.}
The next corollary of the dichotomy completes the picture when finite-energy does not hold under assumption $\textbf{A2}$.
\begin{corollary}
    Let $G_u$ be the directed geodesic graph in direction $u \in \R^2$ constructed in~\citep[Proposition 5.1 and Proposition 5.2]{MR3152744} under assumption \textbf{A2}. Then infinite paths in $G_u$ either coalesce almost surely or $G_u$ contains \biinfinite{} trajectories.
    \label{cor:dichotomy theorem for first and last passage percolation graphs}
\end{corollary}

\begin{proof}
Proposition 5.1 of~\citep{MR3152744} shows that there exists a semi-infinite path from each $x \in \Z^2$. \citep[Proposition $5.2$]{MR3152744} shows that these paths must coalesce when they meet. Together, \citep[Proposition 5.2, Theorem 5.3 and Lemma 6.2]{MR3152744} verify assumptions $1 \AND 2$ of~\Thmref{thm:dichot finish}; in fact, they show that the family of geodesics is asymptotically directed in a sector of angular size at most $\pi/2$. \qed
\end{proof}
    \oldnotes{\citet{MR3535900} use a queuing theory construction to obtain Busemann functions that assumes iid weights.}

The arrows induce a map $T_{\alpha}$ along walks defined by
\begin{equation}
    T_{\alpha}\w = T^{\alpha(\w)}\w  .
\label{eq:translation map along walks}
\end{equation}
The $T_{\alpha}$ map is neither measure preserving nor invertible in general. Along \biinfinite{} trajectories, however, it is both invertible and measure preserving. This observation and~\Thmref{thm:no spec coal} are used in the next theorem. Let $\mathcal{S}$ be the event that the origin is in a \biinfinite{} trajectory. For any $\mathcal{A} \in \mathcal{F}$, let $\Prob_{\alpha}(\mathcal{A}) = \Prob(\mathcal{A} \cap \mathcal{S})$ to obtain the measure space $(\specialpts,\mathcal{F}_{\alpha},\Prob_{\alpha})$.

\begin{definition}[Asymptotic velocity]
   We say that a walk $X \colon \Z^+ \to \Z^d$ has asymptotic velocity if
   \[
       \lim_{k \to \infty} \frac{X(k) \cdot e_i}{k}
   \]
   exists for each $i=1,\cdots,d$.
\end{definition}
\arjunnotes{Indeed, this is the \emph{net} asymptotic velocity in each direction given by
    $$
    \lim_{k \to \infty} \frac{X(w,k) \cdot e_i}{k} = \lim_{k \to \infty} \frac1{k} \sum_{j=1}^k (1_{e_i} - 1_{-e_i})(X(w,j) - X(\w,j-1))
    $$
}
\begin{theorem}\label{thm:measure preserving dynamical system and asym dir}
    The \biinfinite{} trajectories form a measure-preserving $\Z$-system $(\specialpts,\mathcal{F}_{\alpha},\Prob_{\alpha},T_{\alpha})$.
\end{theorem}
Hence, almost surely, ergodic averages converge on all \biinfinite{} trajectories and every \biinfinite{} trajectory has an asymptotic velocity. 
\begin{corollary}
    Suppose the walks in $\mathcal{W}$ satisfy the assumptions in \Thmref{thm:dichot finish}. In $d = 2$, when \biinfinite{} trajectories exist, almost surely, ergodic averages converge on all walks in $\walks{}$ in that configuration, and moreover all walks have the same asymptotic velocity. When the $\Z^2$ system $(\W,\mathcal{F},\Prob,T)$ is ergodic, this direction is deterministic.
    \label{cor:biinfinite trajectories in dim 2}
    \arjunnotes{I wonder if it is enough to say that $\walks$ do not have any loops; it seems like the line crossing assumption should be automatically satisfied since the bi-infinite trajectories have asymptotic direction.}
\end{corollary}
Corollary \ref{cor:biinfinite trajectories in dim 2} is a simple consequence of \Thmref{thm:measure preserving dynamical system and asym dir}, \Thmref{thm:no spec coal}, and \Corref{cor:walks that are not part of biinf}. Note that in general, ergodic averages do not have to have the same limit on all walks.

\begin{oldnote}
In $d>2$ we do not need that the trajectories share an asymptotic direction (almost surely). Indeed consider that for each $(z,0,0)$ we choose that either $(z,a,b)=e_2$ for all $a,b$ or $(z,a,b)=e_3$ for all $a,b$. Note this this $\mathbb{Z}^3$ ergodic but not totally ergodic. 
\begin{remark}
    Suppose we have weights $(\{\w_x\}^{\Z^d},\mathcal{F},\Prob)$ and a stationary set of geodesic walks $\{\walk\}_{z \in \Z^d}$ (formed using the Damron-Hason procedure, say~\citep{MR3152744}). Then, if $\walk_0$ is on a bigeodesic, for any set $A \in \R$, $n^{-1} \sum_{i=1}^n \delta_{\w_{\walk_0(i)}}(A) \to \Prob_{\alpha}(\w_0 \in A)$ as $n \to \infty$ almost surely and in $L^1$. 
\end{remark}
\end{oldnote}
\begin{question}
Under what conditions is the invariant measure in \Thmref{thm:measure preserving dynamical system and asym dir} ergodic? In the context of first- and last-passage percolation, if the walks are Busemann geodesics that form \biinfinite{} trajectories, is there some natural assumption on the weights that ensures that the invariant measure on the \biinfinite{} trajectories is ergodic?
\end{question}

After establishing these nice properties of the bi-infinite trajectories, we show that in the case of almost sure coalescence, none of these properties need to hold. We do this by building an example in $\Z^2$ using a cutting and stacking construction (cutting and stacking was initiated by \citet{MR0247028}).  

\begin{theorem} 
    There exists an ergodic $\mathbb{Z}^2$ dynamical system, $(\Omega, \mathcal{F},\mathbb{P},\{T^z\}_{z \in \Z^2})$, and a stationary arrow map $\alpha \colon \W \times \Z^2 \to \{e_1,e_2\}$ that defines walks from every point on $\Z^2$ such that almost surely, all walks coalesce but no walk has an asymptotic direction.  
    \label{thm:coal but no asym direction in z2}
\end{theorem}

\begin{corollary}
    In the {setting of \Thmref{thm:coal but no asym direction in z2}, the weight function defined in \Eqref{eq:weights in first passage model without asymptotic direction} gives a standard first-passage percolation model with edge-weights on $\Z^2$ and a family of compatible geodesics that coalesce almost surely. Here, with probability $1$, a positive density of points have geodesics with no asymptotic direction.}
     \label{cor:first-passage model with no asymptotic direction or weight distribution}
\end{corollary}
See \Secref{sec:example} for more details about \Corref{cor:first-passage model with no asymptotic direction or weight distribution}. \Remref{rem:no asymptotic weight distribution on geodesics} modifies the weights in \Eqref{eq:weights in first passage model without asymptotic direction} to give a counter-example to C.\,Hoffman's question in the almost surely coalescent setting. Here, the average weight on the geodesic does not converge on a compatible family of geodesics. 

We then use~\Thmref{thm:coal but no asym direction in z2} to construct an example in $\Z^3$ where we have neither almost sure coalescence nor \biinfinite{} trajectories. In other words, the dichotomy theorem no longer holds in $\Z^3$.
\begin{corollary}
    \label{thm:both almost sure coal and biinf walks}
    There exists an ergodic $\mathbb{Z}^3$ system defining walks where in almost every configuration,
    \begin{enumerate}[(i)]
        \item every walk does not have an asymptotic direction,
        \item there are no bi-infinite trajectories, and
        \item there is a positive density of walks that do not coalesce with each other. 
    \end{enumerate}
\end{corollary}

\section{Noncoalescence implies \biinfinite{} trajectories}

Before proceeding with the proof of~\Thmref{thm:dichot finish}, we prove an elementary lemma that we will use repeatedly.
\begin{lemma}
    Let $(\Omega,\mathcal{F},\mathbb{P},\Z^d)$ be an ergodic $\Z^d$ system, and let $M(\w) \subset \Z^d$ be a random, translation covariant $(x \in M(\w) \Leftrightarrow x - z \in M(T^z\w) \, \forall z)$ set of points such that $|M(\w)| \geq 1$ occurs with positive probability. Then, almost surely, $M(\w)$ must have density $\rho > 0$.  
    \label{lem:at least one point implies density}
\end{lemma}
\begin{proof}
    Let $\mathcal{U}_x = \{ \w \colon x \in M(\w) \}$. We claim that $\Prob(\mathcal{U}_x) > 0$. For if not, $\Prob(|M(\w)| > 0) = \Prob( \cup_x \mathcal{U}_x ) = 0$. Then, applying Birkhoff's theorem to the indicator of $\mathcal{U}_\origin$ shows that $M(\w)$ has density equal to $\Prob(\mathcal{U}_\origin)$ almost surely.
\qed
\end{proof}

In a given configuration $\w$, we say that a point $x \in \Z^d$ has a \emph{past of length} $n$ if there is a $z$ such that $\walk_z(\w,n) = x$. Define the random set of points
\[
    P_n(\w) := \{ x \in \Z^2 \colon x \text{ has a past of length $n$} \}.
\]
A point is in a \biinfinite{} trajectory iff it is in $\cap_{n \geq 0} P_n(\w)$. {When $d=2$,} the following proposition shows that the probability that the origin is in a trajectory of length $n$ is bounded below by a constant $\rho > 0$ that is independent of $n$. Therefore $\cap_{n \geq 0} P_n(\w)$ has positive density, and \Thmref{thm:dichot finish} follows.

\begin{proposition}
    \label{prop:past of length n has density rho}
    Under assumptions $1$ and $2$ of \Thmref{thm:dichot finish}, if walks do not coalesce almost surely, there exists a constant $\rho>0$ so that for all $n$, $P_n(\w)$ has density at least $\rho$ almost surely.
\end{proposition}
\Propref{prop:past of length n has density rho} will be proved over several steps in this section.
\begin{openquestion}
    Jon, you removed this old paragraph. Can this be modified to remove the condition that the lines have to be vertical? Given a direction, we choose a countable number of parallel lines in that direction distance one apart. We choose a perpendicular to our direction and identify points in $\mathbb{Z}^2$ with the line this perpendicular direction first intersects. If the point is on a line we identify it with that line. It is straightforward that for each direction there is a well defined density of the point on each line. With this observation, the arguments which we phrase for vertical lines for convenience apply equally well to lines in any direction.
\end{openquestion}

Let $\ver_a^k = \{a\} \times [-k,k]$ and write $\ver_a$ for the entire vertical line with $e_1$ coordinate $a$. Since there is no almost sure coalescence, there must be a pair of points $x,y \in \Z^2$ such that $x$ and $y$ do not coalesce. By assumption 2 in \Thmref{thm:dichot finish}, we may assume without loss of generality that $y$ is of the form $x + k e_2$ for some integer $k \in \Z$. Since the walks cross lines, the walks from {$x$ and $x +k e_2$} must cross $\ver_{x\cdot e_1}$ a last time. The last-crossing points must be of the form $z$, $z + r e_2$ such that $z \cdot e_1 = x \cdot e_1$. So let 
\begin{multline*}
    L_r(\w) := \{ x \in \Z^2 \colon x \AND x + r e_2 \text{ do not coalesce } \\
    \AND \forall k > 0, \walk_{y}(\w,k) \cdot e_1 > x \cdot e_1 , \text{ for } y=x,x+r e_2 \}.
\end{multline*}
\begin{lemma} \label{lem:density of last crossing non coalesing pts} 
    There exists $r$, such that $L_r(\w)$ has density $\xi > 0$ almost surely.
\end{lemma}
\begin{proof}
    This follows directly from \Lemref{lem:at least one point implies density} since $L_r(\w)$ is translation covariant.
\qed
\end{proof}

We say that $x < y$ for two points $x,y \in \ver_a$  if $x \cdot e_2 < y \cdot e_2$. For a walk $\walk_y$, let $M_y \subset \R^2$ be the continuous curve obtained by joining together the vertices in $\walk_y$ by straight line segments.

\begin{lemma} Let $x,y \in \Z^2$ such that $x < y \in \ver_p$ for $p \in \Z$ and both $x$ and $y$ are last crossing points. Then, $\walk_x(k)$ must remain in the closed, unbounded region formed by $\ver_p$ on the left and $M_y$ on the top. If $\walk_x$ touches $M_y$ at some point, it must coalesce with it.
    \label{lem:last points cannot alter vertical ordering in the future}
\end{lemma}
\begin{proof}
Consider the union of $M_y$ and $\ver_p$. It divides the plane into 3 regions: one on the left of $\ver_p(w)$ and two to the right which are divided by the line $M_y$. Now consider $M_x$, the continuous line obtained from the semi-infinite walk $\walk_x$. 
Since $X_x$ remains to the right of $\ver_p$, $M_x \setminus \{x\}$ must intersect the open region below $M_y$ and to the right of $\ver_p$. To exit this region, it must cross $M_y$. If this happens it coalesces with $X_y$. \qed
\end{proof}\arjunhl{intersect the open region below $M_y$ and to the right of $\ver_p$}\arjunhl{But this is not possible since $X_x$ and $X_y$ cannot cross. Therefore, $\walk_x$ must either stay below $\walk_y$ or coalesce with it.}

Consider the set of points $L_r(\w) \cap \ver_p$, and vertically order these points as $\cdots < x_{-1} < x_0 < x_1 < \cdots$ such that $x_0$ is the point with the smallest $e_2$ coordinate in absolute value (with some tie-breaking rule). Let $S_p \subset L_r(\w) \cap \ver_p$ be defined as follows: let $i_0 =0$; for $k \geq 1$, inductively define $x_{i_{k+1}}$ as the point with smallest $y$ coordinate in $L_r \cap \ver_p$ such that $x_{i_{k+1}}- x_{i_k} \geq r$, and analogously define $x_{i_{k}}$ for $k \leq -1$.  Let $S_p := \{ \cdots < x_{i_{-1}} < x_{i_0} < x_{i_1} \cdots \}$. We call $S_p(\w)$ a \emph{separating set} and the points in it \emph{separating points} for the following reason.

\begin{corollary}[of \Lemref{lem:last points cannot alter vertical ordering in the future}]
    \label{cor:stay apart} The walks from two distinct points $x,y \in S_p$ do not coalesce to the right of $\ver_p$.
\end{corollary}
\begin{proof}
    Let $x < y \in S_p$. By \Lemref{lem:last points cannot alter vertical ordering in the future}, $M_x$ must remain in the region bordered by $M_{x + re_2}$ and $\ver_p$. Since $x$ and $x + re_2$ are in $L_r$, $M_x$ must remain in the interior of this region. Similarly, if $y$ is different from $x + re_2$, then $M_y$ must remain in the region bordered by $M_{x + re_2}$ and $\ver_p$. Therefore $x$ and $y$ cannot coalesce.
\qed
\end{proof}

\begin{lemma}\label{lem:plenty} 
    Suppose $\rect_\origin(N)$ has $c (2N+1)^2$ points in $L_r(\w)$. Then, the set $\{ p \colon |S_p \cap \rect_\origin(N)| \geq (c/2r) (2N+1) \}$ has cardinality at least $(c/2) (2N+1)$.
\end{lemma}
\begin{proof}
    Let $M = 2N+1$. Suppose at most $c M/2$ lines have at least $c M/2$ points in $L_r(\w)$, then the number of points in $\rect_\origin(N) \cap L_r(\w)$ is at most
    \[
        \frac{c M}{2} M + \left( 1 - \frac{c}{2} \right) M \frac{c M}{2}
        < c M^2.
    \]
    This contradicts the assumption in the lemma. Therefore, there must be at least $c M/2$ lines with at least $(c/2) M$ points in $L_r(\w)$. On each of these lines, at least $(c/2r) M$ must be in the separating set.
\qed
\end{proof}
\arjunnotes{So you can write this as $x M + (1 - x)M x = f(x)$ for $0 \leq x \leq cM/2$. Since $f'(x) = 2M - 2Mx \geq 0$ for $x \leq 1$, it must take its maximum value at $x = cM/2$}

The following lemma states that each separating point in $\ver^m_k$ corresponds to a unique point on the boundary of $\rect_{{(k,0)}}(n,m+n)$ with past of length at least $n$. This follows from \Corref{cor:stay apart} which says that walks from the separating set must stay apart in the future.
\begin{lemma}Let $S_k(\w)$ be a separating set. For all $m,n > 0$ we have 
    \[
        | P_n(\w) \cap \partial \rect_{(k,0)}(n,m+n) |\geq | S_k(\w)\cap \ver^{m}_{k}(\omega) |.   
    \]
\end{lemma}
\begin{proof}
    {By the line-crossing assumption in \Thmref{thm:dichot finish}}, if $(k,c) \in S_k$ and $ c\in [-m,m]$, then the walk $\walk_{(k,c)}$ must cross $\partial \rect_{(k,0)}(n,m+n)$. Clearly, $\walk_{(k,c)}$ must take at least $n$ steps before it crosses $\partial \rect_{(k,0)}(n,m+n)$, and hence the crossing point must be in $P_n(\w)$. By Corollary \ref{cor:stay apart} distinct points in $S_k$ cross at distinct points in $\partial \rect_{(k,0)}(n,m+n)$.
\qed
\end{proof}
\begin{corollary}Let $\xi$ and $r$ be as in Lemma \ref{lem:density of last crossing non coalesing pts}. For each $n$ and $\epsilon > 0$, there exists $N_0(\e,n)$ so that for all $N > N_0$, 
    $$
    \mathbb{P}\left(\left\{\omega:|P_n(\w) \cap  \partial \rect_{\origin}(N)|> \frac{\xi}{{4}r}N\right\} \right)>(1-\epsilon).
    $$
    \label{cor:with high probability large number of past n points in boundary of rectangle}
\end{corollary}
\begin{proof}
    By the previous lemma it suffices to show that for each $\epsilon>0$ and $n$ there exist $N > n$ and $|k|\leq N-n$ so that  $|S_k(\w) \cap \ver_k^{N-n}|>\frac{\xi}{2r} N$ with probability at least $1-\epsilon$. For $N$ large enough, the ergodic theorem guarantees that there will be $(\xi/2) (2(N-n) + 1)^2$ points in $L_r(\w) \cap \rect_\origin(N-n)$ with probability greater than $1- \e$.   Lemma \ref{lem:plenty} shows that at least one of the $2(N-n) + 1$ vertical lines in $\rect_\origin(N-n)$ must have $(\xi/4r) (2(N-n) + 1)$ separating points. Finally we choose $N_0$ so large such that $2(N-n) > N$ for all $N > N_0$.
\qed
\end{proof}
The next lemma shows that with positive probability, there is a positive density of points ($\geq \xi/256r$) in each configuration that have past of length at least $n$. 
\begin{lemma} Let $\xi$ and $r$ be as in Lemma \ref{lem:density of last crossing non coalesing pts}. There exists $M_0$ so that 
    $$
    \mathbb{P}\left( \left\{ \omega:|P_n(\omega)\cap \rect_{\origin}(M)|> {\frac {\xi} {256 r} (2M+1)^2} \right\} \right)> \frac 1 4
    $$ for all $M>M_0$. 
    \label{lem:positive density of points with positive probability}
\end{lemma}
\begin{proof}
    Let $\mathcal{G}_N:=\{\omega:|P_n(\omega) \cap \partial \rect_\origin(N)|>\frac {\xi}{4r} N\}$ be the event in \Corref{cor:with high probability large number of past n points in boundary of rectangle}. There exists $N$ so that $\Prob(\mathcal{G}_N)>\frac 1 2 $. Therefore by the ergodic theorem for $M > N$, we must have 
    $$
\Prob\left(  \w \colon \left|\{(n,m)\in \rect_\origin(M):T^{(n,m)}(\omega) \in \mathcal{G}_{N}\} \right| \geq \frac{1}{4} (2M+1)^2  \right) \geq \frac14.
$$ 
Hence, for each point $(p,q)\in \rect_\origin(M-N) \cap \mathcal{G}_N$ we have at least $\frac{\xi}{4r} N$ points in $\partial \rect_{(p,q)}(N)$ that are in $P_n(\w)$. Each such point in $P_n(\w)$ can appear in at most $ 8N$ different rectangle boundaries of the form $\partial \rect_{(p,q)}(N)$ for different points $(p,q) \in \rect_\origin(M-N) \cap \mathcal{G}_N$. 
Thus, we obtain the following lower bound on points with past of length $n$ inside a rectangle of size $M$. With probability at least $1/4$,
\begin{align*} 
    | & P_n(\w) \cap \rect_\origin(M) | \\
     & \geq \frac{\left(|\{(n,m) \in  \rect_\origin(M):T^{(n,m)}\omega \in \mathcal{G}_N\}|-|\rect_\origin(M)\setminus \rect_\origin(M-N)|\right)\frac{\xi}{4r} N}{8N} \\
    & \geq \left( (2M+1)^2\frac 1 4-4MN \right) \frac{\xi}{32r} ,
\end{align*}
where the subtraction in the second inequality accounts for boundary effects. By choosing $M$ sufficiently large (given $N$) the lemma follows.
\qed
\end{proof}

\begin{proof}[of \Propref{prop:past of length n has density rho}]
    {Let $\mathcal{P}_n$ be the event that the origin has a past of length $n$. From the ergodic theorem, it follows that if $\Prob(\mathcal{P}_n) \leq \xi/256r - \e$ for small enough $\e > 0$, then 
   \[
       \lim_{M \to \infty} \Prob\left(  \w \colon \frac{1}{(2M+1)^2} | \mathcal{P}_n \cap \rect_\origin(M,M)| > \frac{\xi}{256r} \right) = 0.
   \]
   This contradicts Lemma \ref{lem:positive density of points with positive probability}, and shows that $\Prob(\mathcal{P}_n) \geq \xi/256r$.}
\qed
\end{proof}

Finally, we prove \Corref{cor:walks that are not part of biinf}, which says that walks that are not on \biinfinite{} trajectories must coalesce with \biinfinite{} trajectories. {Recall  $\specialpts := \{ \w \in \W \colon \origin \in \walks, \AND \origin \text{ is \biinfinite{}} \}$}. 

\begin{proof}[of \Corref{cor:walks that are not part of biinf}]
    Let $\walks'(\w) := \left(\walks \cap (\W \setminus \specialpts) \right)(\w)$ be the random set of points with walks that do not coalesce with \biinfinite{} trajectories. These walks inherit stationarity and compatibility from the original walks, and are disjoint from the \biinfinite{} points. If such points exist with positive probability, by \Lemref{lem:at least one point implies density}, $\walks'$ is a positive density subset of $\Z^d$.

   By the ergodic decomposition for the $T^{e_1}$ map, there must be at least one vertical line that has a positive density of points in $\walks'$.  Consider a bi-infinite trajectory such that there is at least one point in $\walks'$ above and below it. This shows that they cannot coalesce, and \Lemref{lem:at least one point implies density} shows that there is a density of such points in $\walks'$.
   
   \arjunnotes{ Jon changed from this: Now, fix some \biinfinite{} trajectory. It can cross each vertical line only a finite number of times. Then, using the ergodic decomposition for the $T^{e_1}$ map, there must be at least one vertical line that has a positive density of points in $\walks'$. Therefore, there must be a pair of points in $\walks'$ in parts of the plane separated by the \biinfinite{} trajectory and hence cannot coalesce.}

    Therefore, \Thmref{thm:dichot finish} says that these walks in $\walks'$ must contain a positive density of \biinfinite{} points. This is a contradiction.
\qed
\end{proof}

\section{Bi-infinite trajectories}
\label{sec: biinfinite trajectories}
In this section we assume that we have bi-infinite trajectories with positive probability. The results in this section apply to all $\Z^d$ systems except for \Corref{cor:biinfinite trajectories in dim 2}, which is only proved for $d=2$.
\begin{definition}
    We say that a point $x \in \Z^d$ is \emph{cataclysmic} if it is a point of coalescence of two distinct \biinfinite{} trajectories in a configuration.
\end{definition}

{The following Lemma is a standard application of the Burton-Keane argument \citep{MR648202}.}
\begin{lemma}Given any bounded rectangle $R \in \Z^d$, the number of points in $\partial R$ crossed by a \biinfinite{} trajectory is at least the number of cataclysmic points in $R$.
\end{lemma}
\begin{proof}
Consider $\walks$ as a directed graph in $\Z^d$ by drawing edges between points $x$ and $y$ in $\Z^d$ if $\walk_x(1) = y$. This graph is a tree since there is exactly one walk from each point in $\walks$, and once walks meet, they coalesce. Let $\mathcal{S}$ denote the set of \biinfinite{} trajectories in $\walks$.  Let $\approx$ be the equivalence relation on $\mathcal{S}$ defined by $\mathfrak{p}$ $\approx$ $\mathfrak{q}$ if $(n,m) \in \mathfrak{p}$ iff $(n,m) \in \mathfrak{q}$ for all $(n,m)\in R$. Let $\mathcal{S}'$ be the set of $\approx$ equivalence classes. Let $\sim$ be the equivalence relation on $\mathcal{S}'$ generated by $\mathfrak{p}\sim \mathfrak{q}$ if there exists a cataclysmic point $(n_1,...,n_d) \in R$ so that $\mathfrak{p}$ coalesces with $\mathfrak{q}$ at $(n_1,...,n_d)$. It suffices to prove the lemma for each $\sim$ equivalence class and so we consider an equivalence class $[\mathfrak{p}] \in \mathcal{S}'/\sim$. Now put an order on the finite set $[\mathfrak{p}]$ and index it by some set $L_{[\mathfrak{p}]}$.  
 After two \biinfinite{} trajectories coalesce, consider the coalesced trajectory to be the largest trajectory that has coalesced. Let $V$ be the set of cataclysmic points in $R$ at which (at least) two distinct elements of $[\mathfrak{p}]$ coalesce. Define $\Phi:V\to L_{[\mathfrak{p}]}$ by 
assigning to each point in $V$, the smallest  trajectory in $[\mathfrak{p}]$ that coalesces with it. By our choice of labeling the continuation of a coalesced trajectory by the largest element that has coalesced, each trajectory in the equivalence class is selected at most once. In other words, $\Phi$ is an injection. Repeating this selection procedure for all equivalence classes, we see that the number of bi-infinite trajectories is an upper bound for the number of cataclysmic points. Each selected trajectory gives at least one distinct crossing of the boundary because the region is bounded and the trajectories have an infinite past. \qed
\end{proof}

\begin{proof}[of~\Thmref{thm:no spec coal}] Assume that at least one cataclysmic point occurs with positive probability. By Lemma \ref{lem:at least one point implies density}, cataclysmic points occur with density $\rho>0$. Choose $N$ so that {$\rho N^d> 2|\partial\rect_\origin(N)|$}. Then with positive probability, there are at least $\frac{\rho}{2} N^d$ cataclysmic points in $\rect_\origin(N)$. By the previous lemma, each cataclysmic point can be associated with a distinct point of $\partial \rect_\origin(N)$. However, there aren't enough points on $\partial \rect_\origin(N)$ to accommodate them all and this is a contradiction. %
\qed
\end{proof}
\begin{remark}
    Using the ergodic decomposition we may relax our assumption to $\mathbb{Z}^d$-preserved measures.
\end{remark}

{Recall the measure space $(\specialpts,\mathcal{F}_\alpha,\Prob_\alpha)$ and the map $T_{\alpha}$ map defined in \eqref{eq:translation map along walks}.
Theorem \ref{thm:measure preserving dynamical system and asym dir} states that $(\specialpts,\mathcal{F}_{\alpha}, \Prob_{\alpha},T_\alpha)$ is a measure-preserving, invertible $\mathbb{Z}$ dynamical system.}
\begin{proof}[of \Thmref{thm:measure preserving dynamical system and asym dir}]
    Since from~\Thmref{thm:no spec coal} the \biinfinite{} trajectories cannot coalesce, $T_{\alpha}$ must be almost surely invertible on the event  $\specialpts$. For each $y \in \directions$, let 
    \[
        \mathcal{V}_y=\{\omega \in \mathcal{S} \colon \alpha(T_\alpha^{-1}\omega)= y \}.
    \]
    Clearly, $\{\mathcal{V}_y\}_{y \in \directions}$ is a partition of $\mathcal{S}$. Thus, for any $\mathcal{E}  \subset \mathcal{S}$,
    \begin{equation*}
        \begin{aligned}
            \Prob_{\alpha}(T_{\alpha}^{-1}(\mathcal{E})) & = \Prob(T_{\alpha}^{-1}(\mathcal{E}) \cap \mathcal{S}) = \Prob(T_{\alpha}^{-1}(\mathcal{E}\cap \mathcal{S}))\\
            & = \sum_{y \in \directions} \Prob(T_{\alpha}^{-1}(\mathcal{E} \cap \mathcal{V}_y)) = \sum_{y \in \directions} \Prob(T^{-y}(\mathcal{E} \cap \mathcal{V}_y)) \\
            & = \sum_{y \in \directions} \Prob(\mathcal{E} \cap \mathcal{V}_y) = \Prob_{\alpha}(\mathcal{E}). 
        \end{aligned}
    \end{equation*}
    The first equality is by the definition of $\Prob_{\alpha}$. The second is because $\specialpts$ is $T_{\alpha}$ invariant. The third follows from the definition of the $\mathcal{V}_y$. The fourth is because on $\mathcal{V}_y$ we have $T^{-1}_{\alpha} = T^{-y}$, and the fifth is because each $T^{-y}$ is $\Prob$ measure preserving. 
\qed
\end{proof}

Let $F = \{ a_1,\ldots,a_m\}$ be any finite set of symbols in $\directions$. Then, we say $F$ appears on a walk if $\{\alpha(T_{\alpha}^{n+i}(\w))\}_{i=1}^m = F$ for some $n \in \Z^+$.
\begin{proposition}
    On each \biinfinite{} trajectory, each finite block appears with some density $\rho \in [0,1]$. 
\end{proposition}
\begin{proof}
    There is an ergodic decomposition of $(\specialpts,\mathcal{F}_\alpha,\Prob_{\alpha},T_\alpha)$. This means that almost every \biinfinite{} trajectory is ergodic for some measure. Therefore, almost surely, each finite block appears with (a possibly $0$) density on all \biinfinite{} trajectories.
\qed
\end{proof}
Applying the previous proposition to each singleton $F_i^\pm = \{\pm e_i\}$, $i=1,\ldots,d$ shows that all \biinfinite{} trajectories have asymptotic velocity. 

Next, we prove \Corref{cor:biinfinite trajectories in dim 2} that shows that in dimension $2$, almost surely, ergodic averages converge on \emph{all walks} in $\walks$, where $\walks$ satisfies the assumptions in \Thmref{thm:dichot finish}. Moreover, all walks in a configuration have the same asymptotic direction. 
\begin{proof}[of \Corref{cor:biinfinite trajectories in dim 2}] 
    \Thmref{thm:measure preserving dynamical system and asym dir} shows that $(\mathcal{S},\mathcal{F}_{\alpha},\Prob_{\alpha},T_\alpha)$ is a measure-preserving dynamical system. In $d=2$, \Corref{cor:walks that are not part of biinf} shows that all walks in $\walks$ must coalesce with \biinfinite{} trajectories, and hence almost surely, ergodic averages converge on all walks.
    
    {Next, we show that all \biinfinite{} walks in a configuration have the same velocity.} We may assume without loss of generality that $\Prob$ is ergodic. If it is not ergodic, we may restrict our attention to an ergodic component of $\Prob$ that has \biinfinite{} trajectories.

    Assume for the sake of contradiction that there is no common asymptotic direction for configurations drawn from $\Prob$. Then, there exists $c$ so that a {$\Prob$}-positive measure set of points are on \biinfinite{} trajectories with slope at least $c$ and a {$\Prob$-positive} measure set of points are on \biinfinite{} trajectories with slope strictly less than $c$. Call these sets $\mathcal{S}_{1}$  and $\mathcal{S}_{2}$ respectively. By the ergodicity of $\Prob$, almost surely, there exist {$i,j,k,\ell \in \Z$} with $i<k$ and $j>\ell$ so that  $T^{(i,j)}\omega \in \mathcal{S}_{2}$ and $T^{(k,\ell)}\omega \in \mathcal{S}_{1}$. Thus, the \biinfinite{} trajectory through $(i,j)$ must coalesce with the \biinfinite{} trajectory through $(k,\ell)$. By Theorem \ref{thm:no spec coal}, this must have zero probability. {Therefore almost surely, all \biinfinite{} trajectories and consequently all walks in $\walks$ have the same asymptotic direction.} If $\Prob$ is ergodic, this direction is deterministic.
\qed
\end{proof}

\section{Examples}\label{sec:example}

\begin{figure}
    \begin{center}
   \def\svgwidth{4.5in}
\begingroup%
  \makeatletter%
  \providecommand\color[2][]{%
    \errmessage{(Inkscape) Color is used for the text in Inkscape, but the package 'color.sty' is not loaded}%
    \renewcommand\color[2][]{}%
  }%
  \providecommand\transparent[1]{%
    \errmessage{(Inkscape) Transparency is used (non-zero) for the text in Inkscape, but the package 'transparent.sty' is not loaded}%
    \renewcommand\transparent[1]{}%
  }%
  \providecommand\rotatebox[2]{#2}%
  \ifx\svgwidth\undefined%
    \setlength{\unitlength}{687.08242911bp}%
    \ifx\svgscale\undefined%
      \relax%
    \else%
      \setlength{\unitlength}{\unitlength * \real{\svgscale}}%
    \fi%
  \else%
    \setlength{\unitlength}{\svgwidth}%
  \fi%
  \global\let\svgwidth\undefined%
  \global\let\svgscale\undefined%
  \makeatother%
  \begin{picture}(1,0.58835524)%
    \put(0,0){\includegraphics[width=\unitlength,page=1]{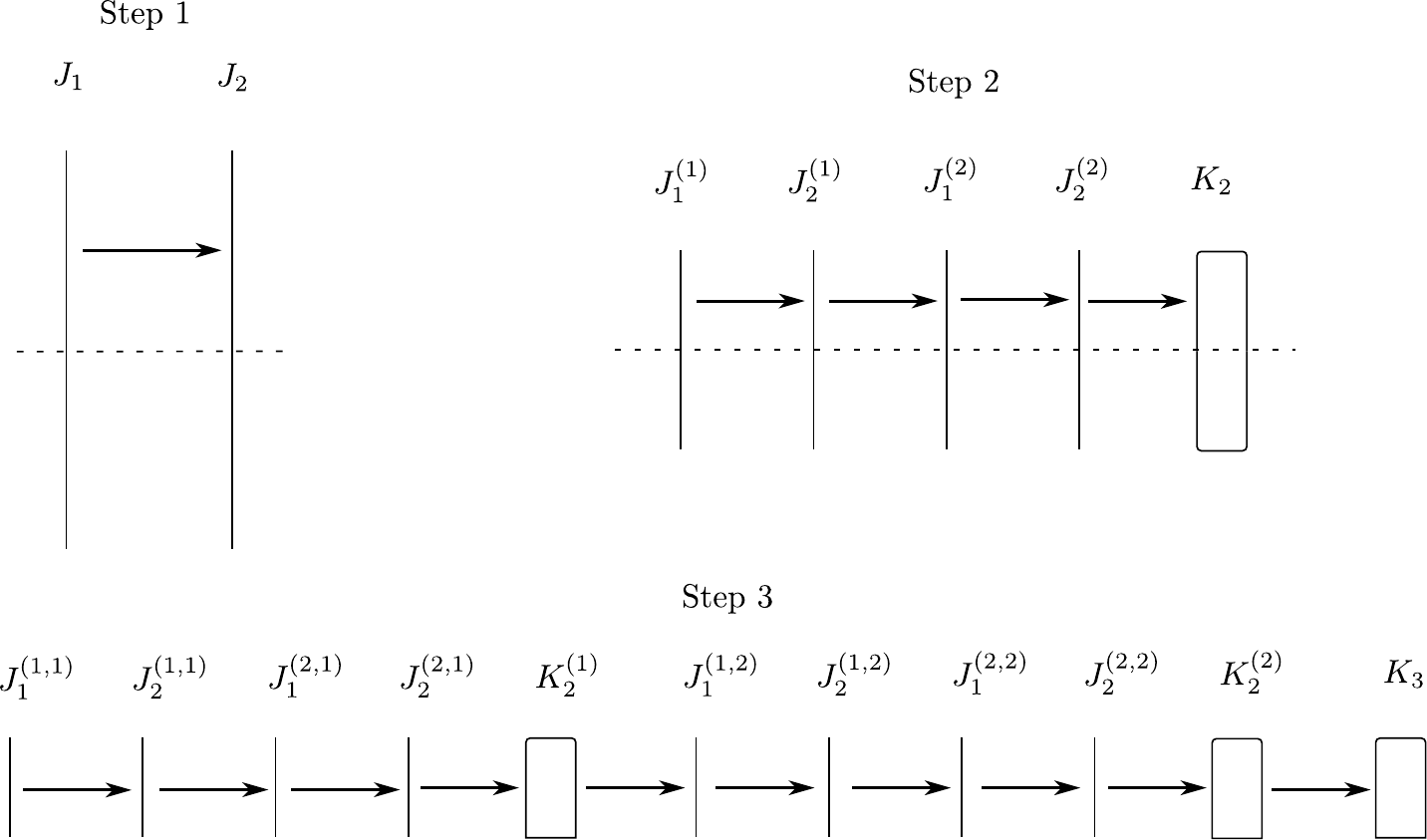}}%
  \end{picture}%
\endgroup%
    \end{center}
   \caption{Three steps of the cutting and stacking construction of the intervals $X$ and $Y$. The figure shows the case $n_1 = n_2 =n_3= 2$. The dotted lines show where the intervals are cut. At each step, the intervals are cut, stacked horizontally, and a new ``spacer'' $K_i$ is appended. The arrows show the mapping $S_i$, $i=1,2$.} 
   \label{fig: cutting and stacking construction}
\end{figure}

We first prove~\Thmref{thm:coal but no asym direction in z2}. We build this example as a product space $X \times Y$, with maps $S_1:X \to X$ and $S_2:Y \to Y$ such that for $(a,b) \in \Z$, we have $T^{(a,b)}(x,y)=(S_1^ax,S_2^by)$ for all $(x,y) \in X \times Y$. $X$ and $Y$ will be intervals in $\R$ with Lebesgue measure.

We build $(S_1,X)$ and $(S_2,Y)$ by cutting and stacking construction. Note that $S_1,S_2$ will be defined as invertible piecewise isometries and therefore Lebesgue measure will be the preserved  measure. In fact they are Rank 1 and therefore ergodic~\citep[Lemma 3]{MR2619397}. 

\textbf{Building $(S_1,X)$:} 
Let $\{ n_i\}_{i=1}^{\infty}$ be a sequence of integers that are at least 2. Let $X=[0,1+\sum_{i=2}^\infty \frac 1 {n_1 \cdots  n_i})$. $X$ is a disjoint union of intervals $(\cup_{i=1}^{n_1} J_i) \cup (\cup_{r=2}^{\infty} K_r)$ that we define below. We consider $n_1$ disjoint intervals of size $\frac 1 {n_1}$, $J_i = [(i-1)/n_1,i/n_1)$ such that $\cup_{i=1}^{n_1} J_i = [0,1)$. 
We begin by defining $S_1$ on $[0,1-\frac 1 {n_1})$. Let $S_1(J_i)=J_{i+1}$ for $i<n_1$. 
 We now define 
$S_1$ on $[0,1)$ so that it agrees with $S_1$ on $[0,1-\frac 1 {n_1})$.   
Subdivide each $J_i$ into $n_2$ intervals of size $\frac{ 1}{n_1n_2}$ called $J_i^{(j)}$  and add an interval $K_2 = [1,1 + 1/n_1 n_2)$. Let 
$S_1(J_i^{(j)})=J_{i+1}^{(j)}$ for all $i<n_1$, $S_1(J_{n_1}^{(j)})=J_1^{(j+1)}$ if $j<n_2$ and lastly $S_2(J_{n_1}^{(n_2)})=K_2$.
 Inductively, we assume that at step $k \geq 2$ we have defined $S_1$ on $[0,1+\sum_{j=2}^{k-1} \frac 1 {n_1\cdot...\cdot n_j})$ which we have divided into intervals $J_i^{(i_2,...,i_k)}, \dots ,K_r^{(j_{r+1},\dots,j_k)},\dots,K_k$ of size $\frac 1 {n_1\cdots n_k}$. An index $i_s$ {or $j_s$} in the superscript runs from $1$ to $n_s$. We now extend the definition of $S_1$ to $[0,1+\sum_{j=2}^{k}\frac 1 {n_1\cdots n_j})$. We add another interval $K_{k+1}$ of size $\frac 1 {n_1\cdots n_{k+1}}$ and subdivide the other intervals into $n_{k+1}$ intervals of size $\frac{1}{n_1\cdots n_{k+1}}$. Call these $J_i^{(i_2,...,i_{k+1})}, \dots, K_r^{(j_{r+1},\dots,j_{k+1})},\dots,K_k^{(\ell)}$, where again, $\ell$ runs from $1$ through $n_{k+1}$. Let
  \begin{equation}
      \begin{array}{ll}
      S_1\,J_i^{(i_2,\ldots,i_{k+1})}=J_{i+1}^{(i_2,\ldots,i_{k+1})} & \text{for } i<n_1, \\
      S_1\,J_{n_1}^{(i_2,\ldots,i_{k+1})}=J_1^{(i_2+1,\ldots, i_{k+1})} & \text{if } i_2<n_2,\\
      S_1\,J_{n_1}^{(n_2,i_3,\ldots,i_{k+1})}=K_2^{(i_3,\ldots,i_{k+1})}  \\
      S_1\,K_r^{(i_{r},\ldots,i_{k+1})}=J_1^{(1,\ldots,i_{r}+1,\ldots,i_{k+1})} & \text{if } i_{r}<n_{r},  \\
      S_1\,K_r^{(n_{r},i_{r+1},\ldots,i_{k+1})}=K_{r+1}^{(i_{r+1},\ldots, i_{k+1})} & \text{otherwise}
      \end{array}
  \end{equation}
  \Figref{fig: cutting and stacking construction} illustrates three steps of the cutting and stacking construction for $n_1 = n_2 = n_3 = 2$.

  \arjunnotes{Suggest we add a proposition here:
      The system is a rank 1 map that is known to be ergodic\cite{MR1436950}.
      Question: Is this system merely measure preserving since we are taking a product of two ergodic dynamical systems. Are the systems weak mixing? Yuri Bakhtin was interested in this. Since this is rank 1, I am guessing its weak mixing but not strong, right?
  }
  \begin{lemma} \label{lem:reach when} For all $x \in X$ {and $r \geq 2$}, $S_1^\ell(x) \in \cup_{i=r}^{\infty} K_i$ for some $0\leq \ell\leq n_1\cdots n_r$. 
  \end{lemma}
 \begin{proof} We prove this by induction on $r$. First we establish the base case of $r=2$.  Observe that if $x \in \cup_{i=2}^\infty K_i$ then it is obviously true. 
      Otherwise, $x\in J_i^{(j)}$ and $S_1^{n_1-i}x\in J_{n_1}^{(j)}$. Now if  $j < n_2$, we have that $S_1^{n_1}(J_{n_1}^{(j)})=J_{n_1}^{(j+1)}$ and if $j=n_2$ then $S_1(J_{n_1}^{(j)})=K_2$. Applying this $n_2-j$ times  we see  $S_1^{(n_2-j)n_1 + 1}(J_{n_1}^{(j)})=K_2$. Combining these, there exists $r\leq n_1-1+(n_2-1)n_1+1 = n_1 n_2$ so that $S_1^r(J_i^{(j)})\subset \cup_{i=2}^\infty K_i$.
 
  The inductive step is similar. Assuming the result for $K_\ell$ we prove it for $K_{\ell+1}$. {If $x \not\in \cup_{i=l}^\infty K_i$, there exists $a\leq n_1\cdots n_\ell - 1$ so that $S_1^{a}(x) \in K_i$ for some $i \geq \ell$.} If $i>\ell$ we are done and so we assume that $S_1^a(x) \in K_\ell^{(j)}$ for some $1 \leq j \leq n_{\ell+1}$. Similar to before $S^{(n_{\ell+1}-j)r + 1 }K_\ell^{(j)}\subset K_{\ell+1}$ where $r\leq n_1 \cdots n_{\ell}$. 
Since $n_1\cdots n_\ell(n_{\ell+1} - 1) + 1 + a\leq n_1 \cdots n_{\ell + 1}$ we have the lemma.   
\qed
\end{proof}
\begin{lemma}\label{lem:reach again} If  $j< n_{i+1}$ then $S_1^{r} K_i^{(j)} \cap \cup_{\ell=i}^\infty K_\ell=\emptyset$ for all $0<r<n_1\cdots n_i$. 
\end{lemma} 
  This is similar to the proof of the previous lemma. 
  
  \vspace{3 mm}
  
  \textbf{Building $(S_2,Y)$:} This is similar. Let $\{m_i\}_{i=1}^{\infty}$ be a sequence of integers that are at least 2. Let $Y=[0,1+\sum_{i=2}^\infty \frac 1 {m_1\cdots m_i}).$
  As before we define intervals $\hat{J}_\ell^{(i_2, \ldots)}$, $\hat{K}_r^{(i_{r+1},\ldots)}$ and the map $S_2$.
  Analogously to before we have the following:
  \begin{lemma}\label{lem:reach when hat} For all $y \in Y$ and $r \geq 2$ we have $S_2^\ell(y) \in \cup_{i=r}^{\infty} \hat{K}_i$ for some $0\leq \ell\leq m_1\cdots m_{r}$.
  \end{lemma}
  \begin{lemma}\label{lem:reach again hat} If $j <m_{i+1}$ then $S_2^{r}\hat{K}_i^{(j)}\cap \cup_{\ell=i}^\infty \hat{K}_\ell=\emptyset$ for all $0<r<m_1\cdots m_i$. 
  \end{lemma}
  For clarity we denote Lebesgue measure on $X$ by $\mu$ and Lebesgue measure on $Y$ by $\hat{\mu}$. Let $\mathcal{B} \AND \hat{\mathcal{B}}$ be the usual Borel $\sigma$-algebras on the intervals $X$ and $Y$.
  \begin{proposition}
    $(X \times Y, \mathcal{B} \times \hat{\mathcal{B}}, \mu \times \hat{\mu}, \{S_1^i \times S_2^j\}_{(i,j) \in \Z^2} )$ is an ergodic $\Z^2$ dynamical system. 
  \end{proposition}
  \begin{proof}
     This is straightforward and included for the reader's convenience. In short, let $A \in \mathcal{B} \times \hat{\mathcal{B}}$ be an invariant set under $S_1^i \times S_2^j$ for all $(i,j) \in \Z^2$. Each section $A_x = \{ y \colon (x,y) \in A \}$ is invariant under $S_2$ and thus $\hat{\mu}(A_x)$ has full or zero measure. But $\hat{\mu}(A_x)$ is measurable function of $x$ that is invariant under $S_1$ and is therefore a constant $\mu$-almost surely. 
  \qed
\end{proof}
  \arjunnotes{Should we include the above proof? I had to prove it for myself.}
  \begin{definition}
      \label{def:arrows map for coalescing walks with no asymptotic direction example}
For any $r > 0$, define $\alpha:X \times Y \to \{e_1,e_2\}$ by 
\begin{equation}
    \alpha(x,y) = \begin{cases}
      e_1 & \text{ if } (x,y) \in \cup_{i=1}^{n_1} J_i \times Y, \text { or } (x,y) \in K_r \times \left(  \cup_{j = r}^\infty  \hat{K}_j \right) \\
        e_2 & \text{ otherwise}
\end{cases}
\label{eq:arrow map definition for cutting and stacking construction}
\end{equation}
This defines an arrow map on $X \times Y$ from which we can construct walks \eqref{eq:direction function definition}.

\end{definition}
\begin{proposition}
    \label{prop: no asymptotic direction}
    If 
    \begin{equation}
        \underset{i \to \infty}{\lim} \frac{n_1 \cdots n_i}{m_1\cdots m_{i-1} i^2} =\infty =\underset{i \to \infty}{\lim} \frac{m_1\cdots m_i}{n_1\cdots n_i i^2}
        \label{eq: condition on mi and ni that ensures that the trajectory fluctuates}
    \end{equation}
    then almost surely, the map $\alpha$ in \eqref{eq:arrow map definition for cutting and stacking construction} defines a trajectory without an asymptotic direction.
\end{proposition} 

\begin{remark}
    The sequences defined by $n_r = 2^{2r -1}$ and $m_r = 2^{2r}$ for $r \geq 1$ satisfy the condition in \Propref{prop: no asymptotic direction}.
\end{remark}

Since we cannot have a $\Z^2$ system with bi-infinite trajectories that do not have asymptotic velocity, \Propref{prop: no asymptotic direction} proves almost-sure coalescence and~\Thmref{thm:coal but no asym direction in z2}. We need a few lemmas to prove~\Propref{prop: no asymptotic direction}.  Let 
\begin{oldnote}
\[
    \hat{G}_r=
    \cup_{j=1}^{r-1} \cup_{\ell_j<m_j }\cup_{\ell_r=1}^{ {m_r}{(1-\frac 1 {r^2})}-1}\cup_{i=1}^{m_1}\hat{J}_i^{(\ell_2,\cdots,\ell_r)} \cup \hat{K}_2^{(\ell_3\dots ,\ell_r)}\cup \dots \cup \hat{K}_{r-1}^{(\ell_r)}.
\]
\end{oldnote}
\begin{equation} \hat{G}_r=Y\setminus  \cup_{j=0}^{m_1\cdots m_r\frac 1 {r^2}} S_2^{-j}( \cup^{\infty}_{l=r} \hat{K}_l).
    \label{eq:set G in the cutting and stacking construction}
\end{equation}
The next lemma says that when $y \in \hat{G}_r$, the walk goes vertically for a period of time that is much longer than its previous horizontal excursion, and hence has a very large vertical fluctuation.
\begin{lemma}\label{lem:travel} If $r\geq 2$ and $y \in\hat{ G}_r$ then {$\walk_\origin(\w,m_1\cdots (\frac{m_r}{r^2}-1))=(a,b)$} where $\w = (x,y)$, and $\frac b a \geq \frac{m_1\cdots ({m_r}\frac 1 {r^2}-1)-n_1\cdots n_r}{n_1\cdots n_r}$.
\end{lemma}
\begin{proof}
We show that the lemma will follow from the fact that if $y \in \hat{G}_r$ then 
\begin{equation} \label{eq:in bad}
S_2^iy \notin \cup^\infty_{\ell=r} \hat{K_\ell}\text{ for } 0\leq i\leq m_1\cdots \left( \frac{m_r}{r^2}-1 \right). 
\end{equation}
By Lemma \ref{lem:reach when}, $S_1^qx \in \cup^{\infty}_{i=r} K_i$ for some $0\leq q\leq n_1\cdots n_r$. 
By Definition \ref{def:arrows map for coalescing walks with no asymptotic direction example}, once $S^q x \in K_r$, subsequent arrows will be in the $e_2$ direction. That is, if $T_{\alpha}^{r}(x,y)$ has its first coordinate in $\cup^{\infty}_{i=r} K_i$, then for $j > r$, $T_\alpha^j(x,y) = T^{e_2}T_\alpha^{j-1}(x,y)$ until the second coordinate of $T_{\alpha}^{j-1}(x,y)$ is in $\cup^\infty_{\ell=r} \hat{K_\ell}$. 
So if \eqref{eq:in bad} holds then $T_{\alpha}^j(x,y)$ moves $e_1$ at most $n_1\cdots n_r$ times in its first $m_1\cdots (\frac{m_r}{r^2}-1)$ steps.

\qed
\end{proof}
Similarly, let
\[ 
    G_r=X\setminus \cup_{j=0}^{n_1\cdot...\cdot n_r\frac 1 {r^2}}S_1^{-j}( \cup_{l=r}^{\infty} K_l )
\]
\begin{lemma}\label{lem:travel 2}  If $r\geq {2}$ and $x \in G_r$ then {$\walk_\origin(\w,n_1\cdots (\frac{n_r}{r^2}-1))=(a,b)$} where $\w = (x,y)$, and $\frac{a}{b}\geq \frac{n_1\cdots (\frac{n_r}{r^2}-1)-m_1\cdots m_{r-1}}{m_1\cdots m_{r-1}}$.
\end{lemma}
\noindent
{The proof of this Lemma is identical to the proof of Lemma \ref{lem:travel}}. 
\begin{lemma}\label{lem:size} If 
 $r \geq 2$, we have $\hat{\mu}(\hat{G}_r^c)\leq \frac2{r^2}$.
\end{lemma}
\begin{proof}
    From \eqref{eq:set G in the cutting and stacking construction} and since we are assuming $m_i\geq 2$ for all $i$, it follows that 
    \begin{align*} 
        \hat{\mu}(\hat{G}_r^c) 
        & \leq \frac{m_1 \ldots m_{r}}{r^2} \sum_{j=r}^{\infty} \hat{\mu}(\hat{K_j}) \\
        & \leq \frac{m_1 \ldots m_{r}}{r^2} \sum_{j=r}^{\infty} \frac1{m_1 \ldots m_{r}(2^{j-r})} \\
        & \leq \frac{2}{r^2}.
    \end{align*}
    \begin{oldnote}
    By construction $\frac{\hat{\mu}(\hat{G}_r)}{\hat{\mu}(Y\setminus \cup_{\ell=r+1}^\infty \hat{K}_\ell)}=\frac{r^2-1}{r^2}-\frac 1 {m_r}$. Since we are assuming $m_i\geq 2$ for all $i$ we have that 
    \[
        \hat{\mu}(Y)= \hat{\mu}(Y\setminus \cup_{\ell=r+1}^\infty \hat{K}_\ell)+\sum_{j=r+1}^\infty \hat{\mu}(\hat{K_j})\leq \hat{\mu}(Y\setminus \cup_{\ell=r+1}^\infty \hat{K}_\ell)+\sum_{j=r+1}^\infty2^{-j}.
    \]
    One more extra step for easy reading would be 
    \[
        \hat{\mu}(Y\setminus \cup_{\ell=r+1}^\infty \hat{K}_\ell)+\sum_{j=r+1}^\infty2^{-j} \leq \hat{\mu}(Y\setminus \cup_{\ell=r+1}^\infty \hat{K}_\ell)+ 2^{-r} \hat{\mu}(Y).
    \]
    \end{oldnote}
\qed
\end{proof}
The next lemma has an analogous proof. 
\begin{lemma}\label{lem:size 2} If 
 $r \geq 2$, we have $\mu(G_r)\leq \frac2{r^2}$. 
\end{lemma}
\begin{proof}[Proof of \Propref{prop: no asymptotic direction}] First observe that by Lemmas \ref{lem:size}, \ref{lem:size 2} and the Borel-Cantelli lemma, 
    $$
    \cup_{i=1}^{\infty}\cap_{r=i}^{\infty} X \times \hat{G}_r  \AND  \cup_{i=1}^{\infty}\cap_{r=i}^{\infty} G_r \times Y
    $$
    have full measure. Next by \eqref{eq: condition on mi and ni that ensures that the trajectory fluctuates} 
    and Lemmas \ref{lem:travel} and \ref{lem:travel 2} we have that any such trajectory approximates both the vertical and the horizontal on infinite sequences of times. 
\qed
\end{proof}


\begin{oldnote}
    I think this is just restating the theorem in~\Secref{sec:main results}. To remove at the end.
\begin{theorem}
    \label{thm:both almost sure coal and biinf walks in the same system}
    There exists an ergodic $\mathbb{Z}^3$ system, $(\Omega', \mathcal{F}',\mathbb{P}')$ and $\alpha$ where almost surely,
    \begin{itemize}
        \item every trajectory does not have an asymptotic direction,
        \item every configuration does not have bi-infinite trajectories, and
        \item we do not have almost sure coalescence. 
    \end{itemize}
\end{theorem}
\end{oldnote}

Proposition~\ref{prop: no asymptotic direction} proves that the trajectories have no asymptotic direction. We turn this into a geodesic walk by assigning weights to {\emph{edges}} on the lattice, and considering the first-passage percolation model on these edge-weights. {Define $w_z \colon X \times Y \to \R$, the weights on the edges $z=e_1 \AND e_2$ by
\begin{equation}
        w_{z}(x,y) =
        \begin{cases}
            \frac 1 2   &    \text{if } \alpha( (x, y) ) = z \\
            1           &    \text{otherwise}
        \end{cases} \quad z=e_1,e_2
        \label{eq:weights in first passage model without asymptotic direction}
    \end{equation}
}
\begin{proof}[Proof of Corollary \ref{cor:first-passage model with no asymptotic direction or weight distribution}] 
    \arjunnotes{If $x \in S_1(\cup_{i=r}^\infty K_i)$ and $y\notin \cup_{i=r+1}^\infty K_i$ let $w(x,y)=1$. Similarly let $w(x,y)=1$ is $y \in S_2(\cup_{i=r}^\infty \hat{K}_i)$ and $x \notin \cup_{i=r}^\infty K_i$. Otherwise let $w(x,y)=\frac 1 2 $.}
{Note that \eqref{eq:weights in first passage model without asymptotic direction} is an everywhere well-defined map.} Our walks in the previous proposition are geodesics for this set of weights because they only cross weights of $\frac 1 2 $. 
\arjunnotes{When we had vertex weights we had to show that there was a non-trivial fraction of vertex weights with value $1$. This showed that the model was nontrivial (not all vertex weights were $1/2$. Now, the following seems unnecessary.
For $(x,y) \in J_1^{(1,2)} \times \hat K_2$, since $S_1^{-1} J_1^{(1,2)} = K_{2}^{(1)}$, we have $\alpha(S_1^{-1}x,y) = e_2$ and $\alpha(x,S_2^{-1}y) = e_1$. This shows that $\w_{\origin}(x,y) = 1$ on $J_1^{(1,2)} \times \hat K_2$, a set of non-zero probability. This means that both weights $1/2$ and $1$ appear with non-zero density in each configuration and hence, the first-passage percolation model is nontrivial.}
\qed
\end{proof}
\begin{remark}[No asymptotic weight distribution on geodesics]
    \label{rem:no asymptotic weight distribution on geodesics}
The previous example can be modified to give an example where there is no asymptotic weight distribution on the geodesic (c.f. \Corref{cor:biinfinite trajectories in dim 2}). 
Let $\hat{w}_{e_1}(x,y) = w_{e_1}(x,y)$, but $\hat{w}_{e_2} = \frac 3 4$ if $w_{e_2}=\frac 1 2$ and $\hat{w}_{e_2} = 1$ otherwise. Walks using these arrows form geodesics since they only cross edges with weight $1/2$ when going horizontally, and edges with weight $3/4$ when going vertically. But horizontal edges weigh $1/2$ or $1$, and vertical edges weigh either $3/4$ or $1$. A calculation similar to the one in \Lemref{lem:travel} and \Lemref{lem:travel 2} shows that the weight on $\walk_\origin(\w,k)$ oscillates between two values $c_1 k$ and $c_2 k$ infinitely often as $k \to \infty$, where $0 < c_1 < c_2$.
\end{remark}

\begin{proof}[Proof of~\Corref{thm:both almost sure coal and biinf walks}] 
        Consider the space $X \times Y \times \Z$ with product measure $(\mu \times \hat{\mu})^{\otimes \Z}$ and translation map $S_1 \times S_2 \times T$ where $T$ is simply the shift on the third coordinate. Let the arrow map be as in \Defref{def:arrows map for coalescing walks with no asymptotic direction example}, where it simply ignores the third coordinate. Observe that almost every trajectory remains in a 2-plane where it behaves as a trajectory from our previous model. Therefore we do not have an asymptotic direction. However, almost every trajectory stays in a 2-plane and so we do not have almost sure coalescence.
\end{proof}

\printbibliography

\end{document}